\newcommand{\xbS}{\Sigma}

\newcommand{\xbs}{\sigma}

\newcommand{\xCd}{\approx}

\newcommand{\xck}{\leq}

\newcommand{\xcp}{\rightarrow}

\newcommand{\xcz}{\Box}

\newcommand{\xDM}{\circ}

\newcommand{\xDd}{\equiv}

\newcommand{\xEI}{\begin{itemize}}
\newcommand{\xEJ}{\end{itemize}}

\newcommand{\xEh}{\begin{enumerate}}
\newcommand{\xEj}{\end{enumerate}}

\newcommand{\xeb}{\prec}
\newcommand{\xec}{\preceq}

\newcommand{\bl}{\begin{lemma} \rm}
\newcommand{\el}{\end{lemma}}
\newcommand{\br}{\begin{remark} \rm}
\newcommand{\er}{\end{remark}}
\newcommand{\be}{\begin{example} \rm}
\newcommand{\ee}{\end{example}}
\newcommand{\bco}{\begin{corollary} \rm}
\newcommand{\eco}{\end{corollary}}
\newcommand{\bc}{\begin{claim} \rm}
\newcommand{\ec}{\end{claim}}
\newcommand{\bfa}{\begin{fact} \rm}
\newcommand{\efa}{\end{fact}}
\newcommand{\bp}{\begin{proposition} \rm}
\newcommand{\ep}{\end{proposition}}
\newcommand{\bd}{\begin{definition} \rm}
\newcommand{\ed}{\end{definition}}
\newcommand{\bcs}{\begin{construction} \rm}
\newcommand{\ecs}{\end{construction}}
\newcommand{\bcd}{\begin{condition} \rm}
\newcommand{\ecd}{\end{condition}}
\newcommand{\bt}{\begin{theorem} \rm}
\newcommand{\et}{\end{theorem}}
\newcommand{\bn}{\begin{notation} \rm}
\newcommand{\en}{\end{notation}}
\newcommand{\bfi}{\begin{bild} \rm}
\newcommand{\efi}{\end{bild}}
\newcommand{\bsta}{\begin{statement} \rm}
\newcommand{\esta}{\end{statement}}
\newcommand{\bcom}{\begin{comment} \rm}
\newcommand{\ecom}{\end{comment}}
\newcommand{\bdia}{\begin{diagram} \rm}
\newcommand{\edia}{\end{diagram}}

\newcommand{\bfc}{\begin{figure}[htb] \begin{center}}
\newcommand{\efc}{\end{center} \end{figure}}

\documentstyle[12pt]{article}
\title{
A REMARK ON UTILITY STREAMS
}

\author{Karl Schlechta
\thanks{
ks@cmi.univ-mrs.fr, karl.schlechta@web.de, http://www.cmi.univ-mrs.fr/ $\sim$ ks
} \\
Laboratoire d'Informatique Fondamentale de Marseille
\thanks{
UMR 6166, CNRS and Universit\'{e} de Provence,
Address: CMI, 39, rue Joliot-Curie, F-13453 Marseille Cedex 13, France
}
}

\date{September 1, 2007}

\begin{document}

\newtheorem{lemma}{Lemma}[section]
\newtheorem{theorem}[lemma]{Theorem}
\newtheorem{proposition}[lemma]{Proposition}
\newtheorem{corollary}[lemma]{Corollary}
\newtheorem{claim}[lemma]{Claim}
\newtheorem{fact}[lemma]{Fact}
\newtheorem{remark}[lemma]{Remark}
\newtheorem{definition}{Definition}[section]
\newtheorem{construction}{Construction}[section]
\newtheorem{condition}{Condition}[section]
\newtheorem{example}{Example}[section]
\newtheorem{notation}{Notation}[section]
\newtheorem{bild}{Figure}[section]
\newtheorem{comment}{Comment}[section]
\newtheorem{statement}{Statement}[section]
\newtheorem{diagram}{Diagram}[section]

\maketitle

\renewcommand{\labelenumi}
  {(\arabic{enumi})}
\renewcommand{\labelenumii}
  {(\arabic{enumi}.\arabic{enumii})}
\renewcommand{\labelenumiii}
  {(\arabic{enumi}.\arabic{enumii}.\arabic{enumiii})}
\renewcommand{\labelenumiv}
  {(\arabic{enumi}.\arabic{enumii}.\arabic{enumiii}.\arabic{enumiv})}

We write down in this very short comment some ideas which occured to the
author during an email discussion with Kaushik Basu on the paper  \cite{BM05}.

The author suggested to consider only finite sequences (which solves the
embedding problem, as the resulting set is countable), and to compare
sequences
of unequal length by repeating them until they have the same length, e.g.
a sequence of length 2 will be repeated 3 times, and a sequence of legth
3 2 times, and the results will then be compared.

Note that the author discussed somewhat related problems in Section 2.2.7
of
 \cite{Sch04}. Considering sums, and not only orders, to evaluate
sequences is
generally difficult, in the sense that often no finite characterizations
are
possible - see again  \cite{Sch04}.

We will write down now a few axioms, which seem reasonable, without
discussion.

We have a domain $X,$ and consider finite, non-empty sequences, noted $
\xbs $ etc.,
with values in $X,$ the set of these sequences will be denoted $ \xbS.$
$X$ has an order
$<,$ $ \xDd $ will express equivalence wrt. this order, and we put
restrictions on a
resulting order $ \xeb $ on $ \xbS,$ with equivalence $ \xCd.$ $ \xck $
and $ \xec $ etc. are defined in the
obvious way.

\bn

$\hspace{0.01em}$


\label{Notation 1.1}

For $ \xbs $ and $ \xbs ' $ of equal length, we write

$ \xbs \xck \xbs ' $ iff all $ \xbs_{i} \xck \xbs '_{i},$

$ \xbs < \xbs ' $ iff $ \xbs \xck \xbs ' $ and for one $i$ $ \xbs_{i}<
\xbs '_{i},$ and finally

$ \xbs << \xbs ' $ iff all $ \xbs_{i}< \xbs '_{i}.$

The double use of $ \xck $ and $<$ will not pose any problem.

$\{x\}$ is the sequence of length 1.

Concatenation is noted $ \xDM.$ For singletons, we may use simple
juxtapposition.

$ \xbs^{n}$ is $ \xbs $ repeated $n$ times.

\en

\paragraph{
Axiom 1.1
}

$\hspace{0.01em}$

(1) Singletons:

(1.1) $x<x' $ $ \xcp $ $\{x\} \xeb \{x' \},$

(1.2) $x \xDd x' $ $ \xcp $ $\{x\} \xCd \{x' \}$

(essentially Pareto).

(2) concatenation:

(2.1) $ \xbs \xDM \xbs \xCd \xbs $ (this expresses essentially that the
mean value is
interesting),

(2.2) $ \xbs ' \xeb \xbs '' $ $ \xcp $ $ \xbs \xDM \xbs ' \xeb \xbs \xDM
\xbs '',$

(2.3) $ \xbs ' \xCd \xbs '' $ $ \xcp $ $ \xbs \xDM \xbs ' \xCd \xbs \xDM
\xbs ''.$

(3) permutation:

$ \xbs \xDM \xbs ' \xCd \xbs ' \xDM \xbs $

(essentially Anonymity).

\bfa

$\hspace{0.01em}$


\label{Fact 1.1}

These axioms allow to deduce:

(4) $ \xbs \xCd \xbs ' $ $ \xcp $ $ \xbs \xCd \xbs \xDM \xbs ' $

(5) $ \xbs \xeb \xbs ' $ $ \xcp $ $ \xbs \xeb \xbs \xDM \xbs ' $

(6) if there are $i,j \xck length( \xbs )=length( \xbs ' ),$ $ \xbs_{i}=
\xbs '_{j},$ $ \xbs_{j}= \xbs '_{i},$ and
$ \xbs_{k}= \xbs '_{k}$ for all other $k,$ then $ \xbs \xCd \xbs ' $ (real
Anonymity)

(7) Weak Pareto:

(7.1) $ \xbs \xck \xbs ' $ $ \xcp $ $ \xbs \xec \xbs ',$

(7.2) $ \xbs < \xbs ' $ $ \xcp $ $ \xbs \xeb \xbs ',$

(7.3) $ \xbs << \xbs ' $ $ \xcp $ $ \xbs \xeb \xbs '.$

and

(8) to compare sequences of different lengths, in the following sense:
When $ \xeb $ and $ \xCd $ are defined between $ \xbs ' $s of equal
length, and Axioms
1-3 hold, then the relation $ \xeb $ (and $ \xCd )$ is determined for
arbitrary sequences.

\efa

\paragraph{
Proof:
}

$\hspace{0.01em}$

Elementary.

(4) by (2.1) and (2.3).

(5) by (2.1) and (2.2).

(6) Let e.g. $ \xbs $ be $ \xbs_{0} \xDM a \xDM \xbs_{1} \xDM b \xDM
\xbs_{3},$
then $ab \xCd ba$ by (3), so $ \xbs_{0} \xDM ab \xCd \xbs_{0} \xDM ba$ by
(2.3) and (3), so
$b \xDM \xbs_{0} \xDM a \xCd a \xDM \xbs_{0} \xDM b$ by (3), so $b \xDM
\xbs_{0} \xDM a \xDM \xbs_{1} \xCd a \xDM \xbs_{0} \xDM b \xDM \xbs_{1}$
by (2.3), etc.

(7) This follows from (1) and repeated use of (2.2) and (2.3).

(8) Let $m:=length( \xbs ),$ $n:=length( \xbs ' ),$ then we obtain by
using (2.1) once,
and (2.3) repeatedly, that $ \xbs^{n} \xCd \xbs $ and $ \xbs '^{m} \xCd
\xbs ',$ but $ \xbs^{n}$ and $ \xbs '^{m}$
have the same length.

$ \xcz $
\\[3ex]

\end{document}